\documentstyle[12pt]{article}
\def\S{{\cal SU}^s_C(r,d)}
\def\ad{\mbox{\rm ad}}
\def\bb{\bf}
\begin{document}
\begin{center}
{\Large \bf  Hecke curves and   Hitchin discriminant}

\bigskip
 {\large \bf Jun-Muk Hwang}
\footnote{Supported  by the Korea Research Foundation Grant
(KRF-2002-070-C00003).} \hspace{3mm} {\large \bf and} \hspace{3mm}
{\large \bf  S. Ramanan}

 \end{center}

\bigskip
\begin{abstract}
Let $C$ be a smooth projective curve of genus $g\geq 4$ over the
complex numbers and ${\cal SU}^s_C(r,d)$ be the moduli space of
stable vector bundles of rank $r$ with a fixed determinant of
degree $d$. In the projectivized cotangent space  at a general
point $E$ of $\S$, there exists a distinguished hypersurface
${\cal S}_E$ consisting of cotangent vectors with singular
spectral curves. In the projectivized tangent space at $E$, there
exists a distinguished subvariety ${\cal C}_E$ consisting of
vectors tangent to Hecke curves in $\S$ through $E$. Our main
result establishes that the hypersurface ${\cal S}_E$ and the
variety ${\cal C}_E$ are dual to each other. As an application of
this duality relation, we prove that any surjective morphism $\S
\rightarrow {\cal SU}^s_{C'}(r,d)$, where $C'$ is another curve of
genus $g$, is biregular. This confirms, for ${\cal SU}^s_C(r,d)$,
the general expectation that a Fano variety of Picard number 1,
excepting the projective space, has no non-trivial self-morphism
and  that morphisms between Fano varieties of Picard number 1 are
rare. The duality relation  also gives simple proofs of the
non-abelian Torelli theorem and the result of Kouvidakis-Pantev on
the automorphisms of $\S$.
\end{abstract}

\section{Introduction}

Any smooth projective variety with Picard group isomorphic to $\bb
Z$ is usually classified into one of three classes, namely {\it
general type, Calbi-Yau or Fano} according as the canonical line
bundle is positive, trivial or negative. Fano varieties are
somewhat special among varieties, and algebraic homogeneous spaces
fall in that class. If we leave out projective spaces, morphisms
between two such varieties of the same dimension seem to be rare
[HM2]. In particular, there is a conjecture, originating from  a
related question of Lazarsfeld,  that there are no nonconstant
self maps of these varieties except automorphisms.

\medskip
Let $C$ be a smooth projective curve of genus $g$ over the complex
numbers and ${\cal SU}^s_C(r,d)$ be the moduli space of stable
vector bundles of rank $r$ with a fixed determinant of degree $d$.
When $r$ and $d$ are coprime, these are smooth Fano varieties with
Picard group $\bb Z$. Thus these provide examples against which
the above kind of conjectures can be tested.

\medskip
 Our main aim in this paper is to prove the following Theorem.

\medskip
{\bf Theorem 5.6} {\it Let $C$ and $C'$ be two smooth projective
curves of genus $g \geq 4$. Let $f: \S \rightarrow {\cal
SU}^s_{C'}(r,d)$ be a surjective morphism. Then $f$ is biregular.
}

\medskip
Note that we do not assume in this theorem that $r$ and $d$ are
coprime, but take only the smooth locus of the varieties in
question. The theorem is perhaps also valid for $g =3$ but our
method does not cover that case. The method synthesizes three
different strands.

\medskip
Firstly, the moduli spaces of vector bundles have been studied by
Hitchin [Hi] from the view-point of symplectic geometry of its
cotangent bundle. On the other hand this study has been used as a
tool to derive results on the moduli spaces themselves in [BNR].
These ideas can be codified in the terms `spectral curves', `Higgs
moduli', `nonabelian theta functions', etc.

Secondly a certain amount of rigidity in the moduli spaces were
established by [NR1] and [NR2] by a study of the geometry of the
moduli spaces. Here the main ingredient is the notion of `Hecke
cycles'. For our purposes it is more fruitful to consider what we
call `Hecke curves' [Hw3].

Finally, the moduli space may be investigated by tools commonly
used in the study of higher-dimensional Fano varieties. This leads
to the study of rational curves on it ([Hw1], [Hw3]) and the Hecke
curves provide the means for doing it. The result quoted above is
obtained by studying interesting relationship between these
aspects.

\medskip
Let us now briefly describe our approach.

Associated to the Hitchin map on the cotangent bundle of $\S$,
there exists a canonically defined hypersurface ${\cal S} \subset
{\bf P}T^*(\S)$ corresponding to twisted endomorphisms of stable
vector bundles whose spectral curves are singular. For a general
point $E \in \S$, the corresponding hypersurface ${\cal S}_E$ in
the projectivized cotangent space ${\bf P}T^*_E(\S)$ will be
called {\it the Hitchin discriminant} at $E$.

On the other hand, there are naturally defined rational curves on
$\S$, which (as referred to above) we call Hecke curves.  For a
general $E \in \S$, let ${\cal C}_E$ be the subvariety of ${\bf
P}T_E(\S)$ consisting of tangent vectors to Hecke curves through
$E$. This subvariety ${\cal C}_E$ will be called the {\it variety
of Hecke tangents} at $E$.

The key point in our proof is the following result which we hope
is sufficiently interesting in itself.

\medskip {\bf Theorem 4.4} {\it Let $g \geq 4$ and let $E$ be a general
point of $\S$. Then the Hitchin discriminant ${\cal S}_E$ is the
dual variety of the variety of Hecke tangents ${\cal C}_E$. }

\medskip
This has other interesting consequences.  It gives simple proofs,
for $g \geq 4$, of non-abelian Torelli theorem (Theorem 5.1) and
the description due to Kouvidakis and Pantev, of the automorphisms
of $\S$ (Theorem 5.4). Our proof of the non-abelian Torelli
theorem is reminiscent of Andreotti's proof of the abelian Torelli
theorem ([An]). Recall that in Andreotti's proof the curve is
recovered as the dual variety of a certain discriminantal
hypersurface associated to the Gauss map of the Riemann theta
divisor. In our proof of non-abelian Torelli theorem, the curve is
recovered from the dual variety of a certain discriminantal
hypersurface associated to the Hitchin map.

\section{Variety of minimal rational tangents}

In this preliminary section, we recall some results concerning
minimal rational curves (cf. [Hw2]). Let $M$ be a smooth
quasi-projective variety of dimension $n$. We will assume that
there exists a component $\cal K$ of the Hilbert scheme of
complete curves on $M$ such that

\smallskip
$(\dagger)$ the subscheme ${\cal K}_y \subset {\cal K}$ consisting
of members of ${\cal K}$ passing through a general point $y\in M$
is a non-empty irreducible smooth projective variety of which
every member is an irreducible smooth rational curve lying in $M$.

\medskip
A member of $\cal K$ is called a {\it minimal rational curve} on
$M$. For a point $y \in M$, let $T_y(M)$ be the tangent space to
$M$ at $y$. Define the {\it tangent morphism}
\begin{eqnarray*} \tau_y: {\cal K}_y &\rightarrow& {\bf P}T_y(M)
\end{eqnarray*} by sending $\ell \in {\cal K}_y$,  a smooth
rational curve $\ell \subset M$,  to
\begin{eqnarray*} \tau_y(\ell) &:=& {\bf P}T_y(\ell). \end{eqnarray*}
For a general member $\ell$ of ${\cal K}_y$, $$T(M)|_{\ell} \cong
{\cal O}(2) \oplus {\cal O}(1)^p \oplus {\cal O}^{n-1-p}$$ where
$p$ is the dimension of ${\cal K}_y$ and ${\cal O}(2)$ corresponds
to $T(\ell)$ ([Hw2, Theorem 1.2]). This implies that $\tau_y$ is
generically finite over its image. The image of $\tau_y$ is
denoted by ${\cal C}_y$ and called the {\it variety of minimal
rational tangents} at the general point $y$ associated to the
family ${\cal K}$.  The following proposition is a consequence of
basic deformation theory.

\medskip
{\bf Proposition 2.1 [Hw2, Theorem 1.4]} {\it Let $\ell$ be a
general member of ${\cal K}_y$ with $$T(M)|_{\ell} \cong {\cal
O}(2) \oplus {\cal O}(1)^p \oplus {\cal O}^{n-1-p}.$$ Then
$\tau_y$ is an immersion at $\ell \in {\cal K}_y$ and the tangent
space to ${\cal C}_y$ at $\tau_y(\ell)$ corresponds to the
subspace of $T_y(M)$ defined by the ${\cal O}(2) \oplus {\cal
O}(1)^p$-part of $T(M)|_{\ell}$.}

\medskip
Recall that when $X$ is an irreducible subvariety of a projective
space ${\bf P}_N$, its dual variety $X^*$ is the irreducible
subvariety of the dual projective space ${\bf P}^*_N$ which is the
closure of the set of hyperplanes containing the projective
tangent space of a smooth point of $X$.  Note that for $\ell$ as
above, $$H^0(\ell, T^*(M)|_{\ell}) \cong H^0(\ell, {\cal O}(-2)
\oplus {\cal O}(-1)^p \oplus {\cal O}^{n-1-p}) = H^0(\ell, {\cal
O}^{n-1-p})$$ are exactly cotangent vectors annihilating ${\cal
O}(2) \oplus {\cal O}(1)^p$-part of $T(M)|_{\ell}$. Also note that
sections of $T^*(M)$ over $\ell$ give smooth rational curves in
$T^*(M)$. As a consequence, we get the following.

\medskip
{\bf Corollary 2.2} {\it Let $\hat{\cal S} \subset T^*(M)$ be the
closure of the union of the smooth rational curves in $T^*(M)$
given by $H^0(\ell, T^*(M))$ as $\ell$ varies over ${\cal K}$. Let
${\cal S} \subset {\bf P}T^*(M)$ be the corresponding projective
subvariety. For a point $y \in M$ let ${\cal S}_y$ be the
intersection ${\cal S} \cap {\bf P}T^*_y(M)$.  Then for general
$y$,  ${\cal S}_y$ is the dual variety of ${\cal C}_y$.}

\medskip
 We recall
the following result from [HM3].

\medskip
{\bf Theorem 2.3 [HM3, Theorem 1]} {\it In the situation above,
the tangent morphism $\tau_y: {\cal K}_y \rightarrow {\cal C}_y$
is birational for a general point $y \in M$.}

\medskip
This was proved in [HM3] when $M$ is a  projective variety, but
the proof there  works even when $M$ is quasi-projective, as long
as the assumption $(\dagger)$ holds.

\medskip
We will also need the following which is essentially [HM1,
Proposition 2].

\medskip
{\bf Proposition 2.4 } {\it Let $M$ and $ \cal K$ be as above.
Suppose there exists  an open subset $A'$ of an abelian variety
$A$ and a generically finite morphism $f:A' \rightarrow M$. Let $y
\in M$ be a general point and $\ell \subset M$ be a general member
of ${\cal K}_y$. Assume that there exists a complete curve $\ell'
\subset A'$ such that $f(\ell') = \ell$. Then the variety of
minimal rational tangents ${\cal C}_y $ is a linear subvariety in
${\bf P}T_y(M)$.}

\medskip
The proof uses the following lemma about curves on abelian
varieties, which is exactly [HM1, Lemma 3].

\medskip
{\bf Lemma 2.5} {\em Let $C_{t}\subset A$ be a $p$-dimensional
irreducible
 family of curves on an $n$-dimensional abelian variety $A$
 passing through a common point $a \in A$.
If  the constructible set in $A$ consisting of the union  of
$C_{t}$'s is of dimension $(p+1)$ and the subspace of
$H^{0}(C_{t}, T^{*}(A))$ consisting of elements annihilating
tangent vectors to $C_{t}$ is of dimension $ \geq n-1-p$ for a
general member $C_{t}$, then the closure of the union of these
curves is a translate of a $(p+1)$-dimensional abelian
subvariety.}

\medskip
{\it Proof of Proposition 2.4}. Let $a \in \ell'$ be a point with
$f(a) = y$.  Note that elements of $H^0(\ell, T^*(M))$ annihilates
the tangent vectors to $\ell$ and $h^0( \ell, T^*(M)) = n-1-p$
where $p$ is the number of ${\cal O}(1)$-factors in
$T(M)|_{\ell}$, or equivalently, the dimension of ${\cal K}_y$.
The pull-back of elements of $H^0(\ell, T^*(M))$ to $H^0(\ell',
T^*(A))$ gives a subspace of dimension $\geq n-1-p$, annihilating
tangent vectors to $\ell'$, because $\ell$ passes through the
general point $y \in M$. By Lemma 2.5, the closure of the union of
all such choices of $\ell'$ is a translate of a
$(p+1)$-dimensional abelian subvariety. In particular, the closure
of  their tangent vectors at $a$ must be a linear subvariety of
${\bf P}T_a(A)$. This implies that ${\cal C}_y$ is a linear
subvariety of ${\bf P}T_y(M)$. $\Box$

\medskip
{\bf Remark 2.6} Since some of our applications, namely, Theorem
5.1 and Theorem 5.3 below,  will be simpler proofs of some results
which have been proved by other means, it is worth pointing out
that the preliminary results reviewed in this section are not so
difficult to prove. The proofs of Proposition 2.1 and Corollary
2.2 are quite straight-forward and use only basic deformation
theory due to Kodaira. Proposition 2.4, whose proof is also easy,
will not be needed for Theorem 5.1 and Theorem 5.3.   The proof of
Theorem 2.3 is more involved, but Theorem 2.3 will be needed in
this paper only when the genus of $C$ is 4.

\section{ Variety of Hecke tangents}

 Let $C$ be a smooth projective curve of
genus $g \geq 4$. Let $\S$ be the moduli space of stable bundles
of rank $r$ with a fixed determinant of degree $d$ over $C$. For
$M= \S$, there exists a family of rational curves satisfying the
condition $(\dagger)$, called Hecke curves. Let us briefly recall
the definition  (see [NR2] and [Hw3] for details).

Let $E \in \S$ be a stable bundle over $C$. Denote by $E^*$ the
dual bundle and ${\bf P}E$ the projectivization consisting of
lines through the origin on each fiber.  For $x \in C$ and $\zeta
\in {\bf P} E^*_x$, consider a new vector bundle $E^{\zeta}$
defined by
\begin{center}$ 0 \longrightarrow E^{\zeta}
\longrightarrow E \longrightarrow (E_x/ \zeta ^{\perp}) \otimes
{\cal O}_x \longrightarrow 0$ \end{center} where $\zeta ^{\perp}$
denotes the hyperplane in $E_x$ annihilated by $\zeta$. Let
$\iota: E^{\zeta}_x \rightarrow E_x$ be the homomorphism between
the fibers at $x$ induced by the sheaf map $E^{\zeta} \rightarrow
E$. The kernel of $\iota$, $Ker(\iota)$, is a 1-dimensional
subspace of the fiber $E^{\zeta}_x$ and its annihilator
$(Ker(\iota))^{\perp}$ is a hyperplane in $(E^{\zeta})^*_x$.  Let
${\bf l}$ be a line in ${\bf P} E^{\zeta}_x$ containing the point
$[Ker(\iota)]$. For each point $l \in {\bf l}$ corresponding to a
1-dimensional subspace $l \subset E^{\zeta}_x$, consider the
vector bundle $\widetilde{E}^l$ defined by
\begin{center} $ 0
\longrightarrow \widetilde{E}^{l} \longrightarrow (E^{\zeta})^*
\longrightarrow [(E^{\zeta})_x^*/l^{\perp}] \otimes {\cal O}_x
\longrightarrow 0$ \end{center} where $l^{\perp} \subset
(E^{\zeta})^*_x$ is the hyperplane annihilating $l$.  This vector
bundle $\widetilde{E}^l$ is stable for each $[l] \in {\bf l}$ if
$E$ is a general point of $\S$ and $g \geq 4$ ([Hw3, Proposition
2]). It is easy to check that for $l= Ker(\iota)$,
\begin{eqnarray*} \widetilde{E}^{Ker(\iota)} &\cong E^*.
\end{eqnarray*} It follows that $\{ (\widetilde{E}^{l})^*; l \in
{\bf l} \}$ defines a rational curve passing through $E$ in $\S$.
A rational curve on $\S$ constructed this way is called a {\it
Hecke curve}. Using [NR2,5.9], one can show that a Hecke curve is
smooth. In view of [NR2,5.16], it is easy to check that a Hecke
curve has degree $2r$ with respect to $K^{-1}_{\S}$.

 On ${\bf P}E^*$, consider the relative cotangent bundle
$\Omega_E$ of the fibration $\varpi: {\bf P}E^* \rightarrow C$.
The projective bundle ${\bf P} \Omega_E$ over ${\bf P}E^*$
 is a smooth projective variety of dimension $2r-2$. The set of
all lines in ${\bf P} E^{\zeta}_x$ containing the point
$[Ker(\iota)]$ is naturally isomorphic to ${\bf P}
(E^{\zeta}_x/Ker(\iota)) \cong {\bf P} \Omega_{E,\zeta}$. In other
words, each point of ${\bf P} \Omega_E$ defines a Hecke curve
through $E$ for a general point $E \in \S$. The argument of
[NR2,5.13] shows that Hecke curves associated to two distinct
points of ${\bf P} \Omega_E$ are distinct rational curves on $\S$.
Thus ${\bf P} \Omega_E$ is naturally isomorphic to the variety of
all Hecke curves through $E$. A simple dimension-counting shows
that Hecke curves are dense in an irreducible component of the
Hilbert scheme of curves on $\S$ ([Hw3, Proposition 3]). It
follows that the component ${\cal K}$ of the Hilbert scheme of
$\S$ corresponding to Hecke curves satisfies the condition
$(\dagger)$, i.e., Hecke curves are minimal rational curves of
$\S$.

Let us describe the tangent morphism associated to Hecke curves
through a general point $E \in \S$.
 Let $\varphi: {\bf P} \Omega_E \rightarrow {\bf P}E^*$ be the
projectivization of $\Omega_E$ and $\xi_E$ be the ${\cal
O}(1)$-bundle of the projectivization so that $\varphi_* \xi_E =
\Omega^*_E$ is the relative tangent bundle of $\varpi$. Recall
that $\varpi_*\Omega^*_E$ is the bundle $ad_E$ of traceless
endomorphisms of $E$.  Let $\pi: {\bf P}\Omega_E \rightarrow C$ be
the composition $\pi= \varpi \circ \varphi$. Note that
\begin{eqnarray*} H^0({\bf P} \Omega_E, \xi_E \otimes \pi^* \omega_C) &=& H^0({\bf P}E^*,
\varphi_*\xi_E \otimes \varpi^* \omega_C) \\&=&
 H^0( {\bf P} E^*, \Omega^*_E \otimes \varpi^*\omega_C) \\
&=& H^0(C, \varpi_* \Omega^*_E \otimes \omega_C) \\ &=& H^0(C,
ad_{E} \otimes \omega_C) \end{eqnarray*} is the dual of the
tangent space of $\S$ at $E$.  Thus the line bundle $\xi_E \otimes
\pi^*\omega_C$ defines a rational map
$$\tau_E: {\bf P} \Omega_E \rightarrow {\bf P}T_E({\cal
SU}^s_C(r,d)).$$
 For a general $E$, this rational map is exactly the tangent morphism
 assigning to each Hecke curve through $E$ its tangent vector at $E$ ([Hw3, Theorem 3]).
  We denote the  image of $\tau_E$ by ${\cal
C}_E$ and call it the {\it variety of Hecke tangents}.

\medskip
{\bf Theorem 3.1} {\it Let $g \geq 5$. Then for a general stable
bundle $E \in {\cal SU}^s_C(r,d)$, the line bundle $\xi_E \otimes
\pi^* \omega_C$ is very ample, i.e., $\tau_E: {\bf P}\Omega_E
\rightarrow {\cal C}_E$ is a biregular morphism.  }

\medskip
{\it Proof}. Write $L$ for $\xi_E \otimes \pi^* \omega_C$. For any
$x \in C$, the line bundle $L$ restricted to the fiber
$\pi^{-1}(x)$ is very ample. Thus $L$ is very ample on ${\bf P}
\Omega_E$ if for any $x , y \in C$, the case $x=y$ included, the
restriction map
\begin{eqnarray*} H^{0}({\bf P}\Omega_E, L) &\longrightarrow& H^0(
\pi^{-1}(x + y), L|_{\pi^{-1}(x+y)}) \end{eqnarray*} is
surjective. {}From the exact sequence \begin{center} $ 0
\longrightarrow L \otimes \pi^* {\cal O}(-x-y) \longrightarrow L
\longrightarrow L|_{\pi^{-1}(x+y)} \longrightarrow 0,$
\end{center} the surjectivity is guaranteed if
\begin{eqnarray*} H^1({\bf P} \Omega_E, L \otimes \pi^* {\cal
O}(-x-y)) &=& H^1( C, ad_E \otimes K_X(-x-y)) \end{eqnarray*} or
its dual $H^0(C, ad_E \otimes {\cal O}(x+y))$ vanishes. Thus
Theorem 3.1 follows from Proposition 3.2 below. $\Box$

\medskip
{\bf Proposition 3.2} {\it Let $\ell$ be a positive integer
satisfying $g \geq \frac{3}{2} \ell +2$. Then for a general stable
bundle $F$ of arbitrary rank and degree $H^0(C, ad_F(D)) =0$ for
any effective divisor $D$ of degree $\ell$.}

\medskip
We need a few lemmas.

{\bf Lemma 3.3} {\it For a general stable bundle $E$ on $C$ of
rank $r$ and degree $d$, $H^0(C, E) = 0$ if $d \leq r(g-1)$.}

\medskip
{\it Proof}. Let us count the dimension of the space of stable
bundles which have non-zero sections.  If $E$ has a non-zero
section, there exists a line subbundle $L \subset E$ with $d':=
\deg(L) \geq 0$. Thus $E$ can be realized as an extension of the
type $$0 \longrightarrow L \longrightarrow E \longrightarrow G
\longrightarrow 0$$ where $L$ is a line bundle of degree $d' \geq
0$ with $H^0(C, L) \neq 0$ and $G$ is a vector bundle of rank
$r-1$ and degree $d^{''}= d-d'$. Since non-stable bundles can be
deformed to stable bundles ([NR1, Proposition 2.6]), we may assume
that $G$ is stable in dimension-counting. Recall that the moduli
space ${\cal U}_C(r,d)$ of semi-stable bundles of rank $r$ and
degree $d$ on $C$ has dimension $r^2(g-1) +1$.  Thus the dimension
of deformation of $G$ is equal to
$$\dim {\cal U}_C(r-1, d^{''}) = (r-1)^2 (g-1) +1.$$ The dimension
of possible choices of the line bundle $L$ is $ \leq d'$. For a
fixed $G$ and a fixed $L$, the dimension of extensions of $G$ by
$L$ is $h^1(C, G^* \otimes L)$. We claim that $H^0(C, G^* \otimes
L) =0$. In fact, assuming that $G= E/L$ for some stable bundle
$E$, if there exists a homomorphism $\eta: G \rightarrow L$, the
composition
$$E  \longrightarrow E/L \stackrel{\eta}{\longrightarrow} L
\longrightarrow E$$ must be identically zero because any
endomorphism of $E$ must be a homothety. Hence $\eta \equiv 0$. It
follows that
\begin{eqnarray*} h^1(C, G^* \otimes L) &=& - \chi(G^* \otimes
L)\\&=&  d^{''} - (r-1) d' + (r-1) (g-1) \\ &=& d -r d' + (r-1)
(g-1). \end{eqnarray*} Thus the space of stable bundles which have
non-zero sections has dimension at most $$d' + (r-1)^2(g-1) +1 +
h^1(C, G^* \otimes L) -1 = (r^2-r) (g-1) + d+ (1-r) d'.$$ Since
$\dim {\cal U}_C(r,d) = r^2(g-1) +1$ and \begin{eqnarray*}
[r^2(g-1) +1]
-[(r^2-r) (g-1) +d+(1-r) d'] &=& r(g-1) -d + (r-1) d' +1 \\
&\geq& r (g-1) -d +1 , \end{eqnarray*}  a general stable bundle
cannot have a non-zero section if $r(g-1) -d \geq 0$. $\Box$

\medskip
{\bf Lemma 3.4} {\it Let $E$ be a general stable bundle of rank
$r-1$ and degree $d$. Assume $\ell$ is a positive integer
satisfying $d
> -(r-1)(g-1-\ell) + \ell$. Then there exists an element $\epsilon
\in H^1(C, E^*)$ such that for any effective divisor $D$ of degree
$\ell$ on $C$, $\psi^D(\epsilon) \neq 0$ where $\psi^D: H^1(C,
E^*) \rightarrow H^1(C, E^*(D))$ is the homomorphism arising from
the short exact sequence
$$ 0 \longrightarrow E^* \longrightarrow E^*(D) \longrightarrow
E^*(D)|_D \longrightarrow 0.$$}

\medskip
{\it Proof}. {}From the exact sequence $$ H^0(C, E^*(D)|_D)
\longrightarrow H^1(C, E^*) \longrightarrow H^1(C, E^*(D)), $$ it
suffices to show that $$h^1(C, E^*) - h^0(C, E^*(D)|_D) - \ell
>0,$$ where $\ell$ is interpreted as the dimension of possible choices of $D$.
Note $h^0(C, E^*) = 0$ by Lemma 3.3 because $\deg(E^*) = -d \leq
(r-1)(g-1). $ Thus $$h^1(C, E^*) = - \chi(E^*) = d +(r-1)(g-1).$$
Also $h^0(C, E^*(D)|_D) = \ell (r-1)$. So $$ h^1(C, E^*) - h^0(C,
E^*(D)|_D) - \ell = d + (r-1)(g-1) - \ell (r-1) -\ell$$ which is
positive from $d > - (r-1)(g-1-\ell) + \ell$. $\Box$

\medskip
{\bf Lemma 3.5} {\it Let $E$ be a general stable bundle on $C$ of
rank $r-1$ and degree $d$.  Let $\ell$ be a positive integer
satisfying $|d| \leq (r-1)(g-1-\ell)$. Assume that $H^0(C,
ad_E(Z)) =0$ for any effective divisor $Z$ of degree $\ell$.
Suppose there exists an extension $F$ of $E$ by ${\cal O},$
 $$ 0 \longrightarrow {\cal O} \longrightarrow F \longrightarrow E
 \longrightarrow 0$$ such that $H^0(C, ad_F(D)) \neq 0$ for some
 effective divisor $D$ of degree $\ell$. Then the extension class
 $[F] \in H^1(C, E^*)$ satisfies $\psi^{D'}([F]) = 0$ for some effective divisor
 $D'$ of length $\ell$.}

 \medskip
 {\it Proof}. Let $\phi:F \rightarrow F(D)$ be a non-zero element
 of $H^0(C, ad_F(D))$. The composition $\beta \circ \phi \circ
 \alpha$ in $$
\begin{array}{ccccccccc} 0 & \longrightarrow & {\cal O} &
\stackrel{\alpha}{\longrightarrow} & F & \longrightarrow & E &
\longrightarrow & 0 \\ & & & & \! \downarrow {\footnotesize \phi}
& & & & \\ 0 & \longrightarrow & {\cal O}(D) & \longrightarrow &
F(D) & \stackrel{\beta}{\longrightarrow} & E(D) & \longrightarrow
& 0 \end{array}$$  defines a section of $E(D)$. Since $\deg(E(D))
= d + \ell (r-1) \leq (r-1)(g-1)$, $E(D)$ cannot have a non-zero
section by Lemma 3.3 and consequently $\beta \circ \phi \circ
\alpha =0$. Thus there exists $\gamma$ satisfying
$$
\begin{array}{ccccccccc} 0 & \longrightarrow & {\cal O} &
\stackrel{\alpha}{\longrightarrow} & F & \longrightarrow & E &
\longrightarrow & 0 \\ & &\! \downarrow {\footnotesize \gamma} & &
\! \downarrow {\footnotesize \phi} & & & & \\ 0 & \longrightarrow
& {\cal O}(D) & \longrightarrow & F(D) &
\stackrel{\beta}{\longrightarrow} & E(D) & \longrightarrow & 0.
\end{array}$$ Let $s \in H^0(C, {\cal O}(D))$ be the section
defined by $\gamma$. Consider $$\phi':= \phi - I_F \cdot s \;\; :
F \rightarrow F(D)$$ where $I_F$ denotes the identity map of $F$.
Then $\phi'$ annihilates the subbundle $\alpha: {\cal O} \subset
F$, inducing a non-zero homomorphism $\zeta: E \rightarrow F(D)$
satisfying
$$
\begin{array}{ccccccccc} 0 & \longrightarrow & {\cal O} &
\stackrel{\alpha}{\longrightarrow} & F & \longrightarrow & E &
\longrightarrow & 0 \\ & & \! \downarrow {\footnotesize 0} & & \!
\downarrow {\footnotesize \phi'} & & \;\;\; \downarrow {\footnotesize \beta \circ \zeta} & & \\
0 & \longrightarrow & {\cal O}(D) & \stackrel{\alpha'}
\longrightarrow & F(D) & \stackrel{\beta}{\longrightarrow} & E(D)
& \longrightarrow & 0.
\end{array}$$
If $\beta \circ \zeta =0$, then there exists $\xi: E \rightarrow
{\cal O}(D)$ such that $\alpha' \circ \xi = \zeta$. Then we get a
non-zero element $\xi^*$ in $H^0(C, E^*(D))$ which is not possible
by Lemma 3.3 because $$\deg(E^*(D)) = -d + \ell (r-1) \leq
(r-1)(g-1).$$ Thus $\beta \circ \zeta \neq 0$. By the assumption
$H^0(C, ad_E (D))=0$, we conclude that $$ \beta \circ \zeta = I_E
\cdot s'$$ for some non-zero $s' \in H^0(C, {\cal O}(D))$. Let
$D'$ be the effective divisor defined by $s'$. We claim that
$\psi^{D'}([F]) =0$ in $H^1(C, E^*(D'))$, which proves the lemma.

To prove the claim, let us recall the definition of the extension
class $[F]$ and $\psi^{D'}([F])$. Let $\delta: H^0(C, E^*\otimes
F) \rightarrow H^1(C,E^*)$ be the boundary map associated to the
short exact sequence $$0 \longrightarrow E^* \longrightarrow E^*
\otimes F \longrightarrow E^* \otimes E \longrightarrow 0.$$ Then
$[F] := \delta(I_E)$ for the identity map $I_E \in H^0(C,
E^*\otimes E)$. The multiplications by $s'$
$$
\begin{array}{ccccccccc} 0 & \longrightarrow & E^* &
\longrightarrow & E^* \otimes F & \longrightarrow & E^* \otimes E
& \longrightarrow & 0 \\ & & \; \! \downarrow {\footnotesize \cdot
s'} & & \; \! \downarrow {\footnotesize \cdot s'} & & \; \!
\downarrow {\footnotesize
\cdot s'} & & \\
0 & \longrightarrow & E^*(D') & \longrightarrow & E^* \otimes
F(D') & \longrightarrow & E^* \otimes E(D') & \longrightarrow & 0.
\end{array}$$ induces a commutative diagram
$$ \begin{array}{ccccc} H^0(C, E^* \otimes F) & \longrightarrow & H^0(C, E^* \otimes E) &
\stackrel{\delta}{\longrightarrow} & H^1(C, E^*) \\ \! \downarrow
{\footnotesize \cdot s'} & & \! \downarrow {\footnotesize \cdot
s'} & & \! \downarrow {\footnotesize \psi^{D'}}
\\H^0(C, E^* \otimes F(D')) & \stackrel{\tilde{\beta}}{\longrightarrow}
&  H^0(C, E^* \otimes E(D')) &
\stackrel{\delta^{s'}}{\longrightarrow} & H^1(C, E^*(D'))
\end{array}. $$ It follows that $$\psi^{D'}([F]) =
\psi^{D'}(\delta (I_E)) = \delta^{s'}(I_E \cdot s').$$ But we know
that $$I_E \cdot s' = \beta \circ \zeta = \tilde{\beta}(\zeta)$$
for some $\zeta \in H^0(C, E^* \otimes F(D))$. Thus
$\psi^{D'}([F]) =0$. $\Box$

\medskip
{\bf Lemma 3.6} {\it Let $ r\geq 2,  \ell \geq 1$ and $d$  be
integers satisfying $$ -(r-1)(g-1-\ell) + \ell < d \leq (r-1)
(g-1-\ell).$$ Suppose for a general stable bundle $E$ of rank
$r-1$ and degree $d$, $H^0(C, ad_E (D)) =0$ for any effective
divisor $D$ of degree $\ell$. Then for a general stable bundle $F$
of rank $r$ and degree $d$, $H^0(C, ad_F(D)) =0$ for any effective
divisor $D$ of degree $\ell$.}

\medskip
{\it Proof}. By Lemma 3.4, we can choose $[F_0] \in H^1(C, E^*)$
such that $\psi^{D}([F_0]) \neq 0$ for any effective divisor $D$
of degree $\ell$. Then by Lemma 3.5, $H^0(C, ad_{F_0}(D)) =0$ for
any effective divisor $D$ of degree $\ell$. By [NR1, Proposition
2.6], $F_0$ can be approximated by stable bundles, i.e.,  there
exists a flat family of bundles $\{F_t, t \in T \}$ parametrized
by an affine curve $T$ with a base point $0\in T$ such that $F_t$
is stable for $t \neq 0$.  On $C^{(\ell)} \times T $ where
$C^{(\ell)}$ is the set of effective divisors of degree $\ell$,
consider the loci of points  $(D, t) \in C^{(\ell)} \times T$ such
that $H^0(C, ad_{F_t}(D)) \neq 0$. This loci is a closed
subvariety of $C^{(\ell)} \times T$ and is disjoint from
$C^{(\ell)} \times \{0\}$ since $H^0(C, ad_{F_0}(D)) =0$ for all
$D \in C^{(\ell)}$.  It follows that there exists $t \neq 0$ such
that $H^0(C, ad_{F_t}(D)) =0$ for all $D \in C^{(\ell)}$. In
particular, for a general stable bundle $F$, $H^0(C, ad_{F}(D)) =
0$ for all $D \in C^{(\ell)}$. $\Box$

\medskip
{\it Proof of Proposition 3.2}. The proof is by induction on the
rank $r$ of $F$. If $r=1$ this is obvious. Assume that the result
holds for a general stable bundle $E$ of rank $r-1$ and degree $d
= \deg(F)$. By Lemma 3.6, the result follows if $$
-(r-1)(g-1-\ell) + \ell <d \leq (r-1)(g-1-\ell).$$ Note that there
are $2(r-1)(g-1-\ell) -\ell $ consecutive  integers $d$ satisfying
this and
$$2(r-1)(g-1-\ell) -\ell  \geq r \;\; \mbox{ for } r \geq 2 \mbox{
and } g \geq \frac{3}{2} \ell +2.$$ If $H^0(C, ad_F(D)) =0$ for
some vector bundle $F$ then $H^0(C, ad_{F'}(D))=0$ for any vector
bundle $F'$ of the form $F'=F \otimes L$ for a line bundle $L$.
Thus  we may assume that the degree $d$ of $F$ belongs to any set
of $r$ consecutive integers. This finishes the proof of
Proposition 3.2. $\Box$

\medskip
{\bf Theorem 3.7} {\it Let $g =4$. Then for a general stable
bundle $E \in {\cal SU}^s_C(r,d)$,  $\tau_E:{\bf P}\Omega_E
\rightarrow {\cal C}_E$ is a birational morphism and is unramified
in a neighborhood of a general fiber of  $ \pi: {\bf P}\Omega_E
\rightarrow C$. }

\medskip
{\it Proof}  The birationality of $\tau_E$ over its image is from
Theorem 2.3. That $\tau_E$ is unramified in a neighborhood of a
general fiber of $\pi$ follows from Proposition 3.8 below, in the
same way that Theorem 3.1 followed from Proposition 3.2. $\Box$

\medskip
{\bf Proposition 3.8} {\it Let $g \geq 4$ and $F$ be a general
stable bundle of arbitrary rank and degree. Then there exists a
point $x \in C$ such that $H^0(C, ad_F(2x)) =0$.}

\medskip For the proof of Proposition 3.8, we need the following
three lemmas, Lemmas 3.9, 3.10 and 3.11, which are just slight
modifications of Lemmas 3.4, 3.5 and 3.6, respectively.

\medskip
{\bf Lemma 3.9} {\it Let $E$ be a general stable bundle of rank
$r-1$ and degree $d$ satisfying $ d > -(r-1)(g-3) $. Then there
exists an element $\epsilon \in H^1(C, E^*)$ such that for a given
$x \in C$, $\psi^{2x}(\epsilon) \neq 0$ where $\psi^{2x}: H^1(C,
E^*) \rightarrow H^1(C, E^*(2x))$ is as defined in Lemma 3.4 with
$D = 2x$. }

\medskip
{\it Proof}. As in the proof of Lemma 3.4, it suffices to show
$$h^1(C, E^*) - h^0(C, E^*(2x)|_{2x})  >0. $$ But this is obvious
from  $$h^1(C, E^*) = d + (r-1)(g-1), \;\; h^0(C, E^*(2x)|_{2x}) =
2(r-1), $$ as  in the proof of  Lemma 3.4. $\Box$

\medskip
{\bf Lemma 3.10} {\it Let $x$ be a point satisfying $h^0(C, {\cal
O}(2x)) = 1$, which is certainly true for a general $x \in C$. Let
$E$ be a vector bundle of rank $r-1$ and degree $d$ satisfying
$|d| \leq (r-1)(g-3)$. Assume that $H^0(C, ad_E(2x)) =0$. Suppose
 $F$ is an extension of $E $ by ${\cal O}$ with $H^0(C, ad_F(2x)) \neq 0$.
Then the extension class $[F] \in H^1(C, E^*)$ satisfies
$\psi^{2x}([F]) =0$.}

\medskip
{\it Proof}. The proof of Lemma 3.5 works almost verbatim. It
suffices to replace the divisors $D$ and $D'$ by $2x$ and the
sections $s$ and $s'$ by the unique section (up to scalar) of
${\cal O}(2x)$. $\Box$

\medskip
{\bf Lemma 3.11} {\it Let $r \geq 2$ and $d$ be integers
satisfying
$$ -(r-1)(g-3)  < d \leq (r-1)(g-3).$$ Suppose for a general
stable bundle $E$ of rank $r-1$ and  degree $d$, $H^0(C, ad_E(2x))
=0$ for some $x \in C$. Then for a general stable bundle $F$ of
rank $r$ and $\det(F) = \det (E)$, $H^0(C, ad_F(2x)) =0$. }

\medskip
{\it Proof}. A simple modification of the proof of Lemma 3.6
works. It suffices to take $\{ F_t \}$ with $\det(F_t) = \det
(E)$, replace $C^{(\ell)}$ by $C$ and use Lemma 3.9 and Lemma 3.10
in place of Lemma 3.4 and Lemma 3.5, respectively. $\Box$

\medskip
{\it Proof of Proposition 3.8}. The proof is by induction on the
rank $r$ of $F$ as in the proof of Proposition 3.2. If $r=1$, it
is obvious. Assume that the result holds for a general stable
bundle $E$ of rank $r-1$ and  $\det(E) = \det(F)$. By Lemma 3.11,
the result follows if $$ -(r-1)(g-3)  < d \leq (r-1)(g-3).$$ Note
that there are $2(r-1)(g-3) $ consecutive integers $d$ satisfying
these inequalities and $$ 2(r-1)(g-3) \geq r  \;\; \mbox{ if }
\;\; r \geq 2 \; \mbox{ and } \; g \geq 4 .$$ If a vector bundle
$F$ satisfies $H^0(C, ad_F(2x)) =0$ for some $x \in C$, then
$H^0(C, ad_{F'}(2x)) =0$ for any vector bundle $F'$ of the form
$F' =F \otimes L$ for a line bundle $L$. Thus we may assume that
the degree $d$ of $F$ belongs to any set of $r$ consecutive
integers. This finishes the proof of Proposition 3.8. $\Box$

\section{Hitchin discriminant and its dual variety}

Let us briefly recall the definition of the Hitchin map and
spectral curves. See [BNR], [Hi] and [KP] for details. As before,
$C$ is a smooth projective curve of genus $\geq 4$.
  Let $$W:= H^0(C, \omega_C^{\otimes 2}) \oplus \cdots \oplus H^0(C,
  \omega_C^{\otimes r})$$ be the space of characteristic
  polynomials and $h: T^*(\S) \rightarrow W$ be the {\it Hitchin
  map}
  defined by $$h(\theta) := (s_2(\theta), \ldots, s_r(\theta))$$
  where for $\theta \in T^*_E(\S) = H^0(C, ad_{E} \otimes \omega_C),$
  $$ s_i(\theta) := (-1)^{i} \mbox{tr}(\wedge^i \theta).$$
Let $K_C$ be the total space of the canonical line bundle
$\omega_C$ and $\alpha: K_C \rightarrow C$ be the natural
projection. For an element $s =(s_2, \ldots, s_r) \in W$, the {\it
spectral curve} $C_s$ associated to $s$ is the curve in the total
space $K_C$ defined by the equation
$$x^r + s_2 x^{r-2} + \cdots + s_{r-1} x + s_r $$ where $x$ is the
tautological section of $\alpha^* \omega_C$.   Let ${\cal D}
\subset W$ be the set of characteristic polynomials with singular
spectral curves.  The following two facts are standard.

\medskip
{\bf Proposition 4.1 [KP, Corollary 1.5 and Remark 1.7]}  {\it
${\cal D}$ is an irreducible hypersurface in $W$ and for a general
point $s \in {\cal D}$, $C_s$ is an integral curve with a unique
ordinary double point over a general point of $C$. }

\medskip
{\bf Proposition 4.2 [BNR, 3.6 and 3.7]} {\it If $s \in W$ has an
integral spectral curve, then $h^{-1}(s)$ is irreducible and  for
a general $\alpha \in h^{-1}(s)$ regarded as an element of $
H^0(C, ad_E \otimes \omega_C)$ for some $E \in \S$, each
eigenvalue of $\alpha_x: E_x \rightarrow E_x \otimes \omega_C$ has
one-dimensional eigenspace for each $x \in C$.  If furthermore the
spectral curve is smooth, i.e.,
 $s \in W\setminus{\cal D}$,  then $h^{-1}(s)$ is an open subset of an abelian variety and
is dominant over $\S$. }

\medskip
{\bf Proposition 4.3} {\it The hypersurface $h^{-1}({\cal D})$ in
$T^*(\S)$ is irreducible.}

\medskip
{\it Proof}. Since $h^{-1}(s)$ for a general $s \in {\cal D}$ is
irreducible by Proposition 4.1 and Proposition 4.2, there exists a
unique irreducible component $S_1$ of $h^{-1}({\cal D})$ which is
dominant over ${\cal D}$. Suppose there exists another component
$S_2$ which is not dominant over ${\cal D}$. We will get a
contradiction.

For $E \in \S$, let us denote the restriction of $h$ to the
cotangent space $T^*_E(\S)$ by  $$h_E: T^*_E(\S) \longrightarrow
W.$$ There is a natural ${\bb C}^{\times}$-action on $T^*_E(\S)$
by the scalar multiplication and a natural ${\bb
C}^{\times}$-action on $W$ by the weighted scalar multiplication.
Clearly, $h_E$ is equivariant with respect to these actions of
${\bb C}^{\times}$. Suppose that $T^*_E(\S) \cap h^{-1}(0) = 0$.
 Then $h_E$ descends to a morphism
$$\check{h}_E : {\bf P} T^*_E(\S) \longrightarrow {\bf P}_{\rm
weight} W$$ where ${\bf P}_{\rm weight}W$ is the weighted
projective space obtained as the quotient of $W \setminus 0$ by
the weighted ${\bb C}^{\times}$-action. This $\check{h}_E$ must be
a finite morphism. It follows that $h_E$ is a finite morphism.

If $S_2$ intersects $T_E^*(\S)$ for some $E$ with $T^*_E(\S) \cap
h^{-1}(0)=0$, the intersection $S'_2:= S_2 \cap T^*_E(\S)$ must be
a hypersurface in $T^*_E(\S)$. But then $S'_2$ is dominant over
the irreducible hypersurface ${\cal D}$ because $h$ is finite on
$T^*_E(\S)$, a contradiction. Thus the image $pr(S_2)$ under the
natural projection $pr: T^*(\S) \rightarrow \S$ is contained in
the subvariety
$${\cal N} : = \{ E \in \S:  \dim(T_E^*(\S) \cap h^{-1}(0)) \geq 1 \}.$$
Recall that ${\cal N} \neq \S$ by [La]. Thus $pr(S_2)$ is a
hypersurface in $\S$ and $S_2 = pr^{-1}(pr(S_2))$. Since
$\dim(h^{-1}(0)) = \dim (\S)$ from [La],
$$\dim( T^*_E(\S) \cap h^{-1}(0)) =1$$ for a general $E \in pr(S_2)$.
 Thus $h_E: T^*_E(\S)
\rightarrow W$ must have general fiber dimension $\leq 1$. This
implies that $h_E(T^*_E(\S))$ is a hypersurface in $W$. Since
$T^*_E(\S) \subset S_2$, this is a contradiction to the fact that
$S_2$ is not dominant over ${\cal D}$.  $\Box$

\medskip
 Let ${\cal S}$ be the
hypersurface in ${\bf P} T^*(\S)$ corresponding to  $h^{-1}({\cal
D})$ in $T^*(\S)$. For a general point $E \in \S$, the
hypersurface ${\cal S}_E := {\cal S} \cap {\bf P}T^*_E(\S)$ will
be called the {\it Hitchin discriminant} at $E$.

\medskip
Recall that when $X\subset {\bf P}_N$ is a smooth subvariety, its
dual variety is the subvariety of the dual projective space ${\bf
P}_N^*$ corresponding to singular hyperplane sections of $X$.
Suppose the normalizaiton $\hat{X}$ of $X$ is smooth and
$\tau:\hat{X} \rightarrow X \subset {\bf P}_N$ is the
normalization morphism. Then  $X^*$ is the  closure of the set of
hyperplanes containing the projective tangent space of a point of
$X$ where $\tau$ is an immersion. This observation will be used
implicitly in the proof of the following theorem, for the case of
$g=4$.

\medskip
{\bf Theorem 4.4} {\it Assume $g \geq 4$. Let $E \in \S$ be a
general point. Then the Hitchin discriminant ${\cal S}_E \subset
{\bf P}T^*_E(\S)$ is the dual variety of the variety of Hecke
tangents ${\cal C}_E \subset {\bf P} T_E(\S)$. In other words,
${\cal S}$ defined above agrees with ${\cal S}$ in Corollary 2.2.}

\medskip
{\it Proof}.  Let $\theta \in h^{-1}({\cal D})$ be a general
point. Then $\theta: E \rightarrow E \otimes \omega_C$ is an
endomorphism of a general stable bundle $E$ such that its spectral
curve $C_{h(\theta)}$ has a unique ordinary double point
singularity which lies over a general point of $C$. It suffices to
show that $\theta \in {\bf P}T^*_E(\S)$ belongs to the dual
variety of ${\cal C}_E$. By Proposition 4.2, for each $x \in C$,
each eigenvalue of $\theta_x$
 has a 1-dimensional eigenspace. Thus we have a curve
$$C_{\theta} \subset {\bf P}E^*$$ biregular to the spectral curve
$C_{h(\theta)}$ corresponding to the 1-dimensional eigenspaces.

Let $\Omega^*_E$ be the relative tangent bundle of the projective
bundle $\varpi: {\bf P}E^* \rightarrow C$. Recall that when an
endomorphism of a vector space $V$ is regarded as a vector field
on ${\bf P}V$, the zero set of the vector field corresponds to the
set of eigenvectors of the endomorphism.  When $\theta$ is
regarded as a vertical vector field on ${\bf P}E^*$ twisted by
$\varpi^* \omega_C$ via the isomorphism
$$H^0(C, \ad_E \otimes \omega_C) \cong H^0({\bf P}E^*, \Omega_E^*
\otimes \varpi^* \omega_C),$$  $\theta$ vanishes exactly on
$C_{\theta}$. Thus, when we regard it as a section of $\xi_E
\otimes \pi^* \omega_C$ on ${\bf P} \Omega_E$, it defines an
element of the linear system $|\xi_E \otimes \pi^* \omega_C|$ with
a singular point lying over the singular point of $C_{\theta}$ by
Lemma 4.5 below. This implies that $\theta$ belongs to the  dual
variety of ${\cal C}_E$ because $\tau_E$ is an immersion over a
general point of $C$ by Theorem 3.1 and Theorem 3.7. $\Box$

\medskip
{\bf Lemma 4.5} {\it Let $$ V= a_1(z_1, \ldots, z_n)
\frac{\partial}{\partial z_1} + \cdots + a_{n-1}(z_1, \ldots, z_n)
\frac{\partial}{\partial z_{n-1}}$$ be a  holomorphic vector field
on the polydisc ${\Delta}^n$. Assume that the zero set of the
vector field
$$a_1(z) = \cdots = a_{n-1}(z) =0$$ is a curve with a singularity at $0$.
Let $\Omega$ be the relative cotangent bundle of the projection
$p:\Delta^n \rightarrow \Delta, p(z_1, \ldots, z_n) = z_n$.
 Then the hypersurface in ${\bf P}(\Omega) \cong {\bf P}^{n-2} \times \Delta^n$
 defined by  $$ a_1(z_1, \ldots, z_n)y_1 + \cdots + a_{n-1}(z_1, \ldots, z_n)
 y_{n-1} =0 $$ for the homogeneous coordinates $[y_1: \cdots: y_{n-1}]\in {\bf P}^{n-2}$
 has a singular point over $z_1=\cdots=z_n=0$.}

 \medskip
 {\it Proof}.
Since $$a_1(z) = \cdots = a_{n-1}(z) =0$$ is a curve with a
singular point at $0\in \Delta^n$, the matrix $(\frac{\partial
a_i}{\partial z_j})|_{z=0}$ has rank $\leq n-2$, by Jacobian
criterion of smoothness. Thus  there exist complex numbers $c_1,
\ldots, c_{n-1},$ with $c_i \neq 0$ for some $i$,  satisfying
$$\sum_{i=1}^{n-1} \left( \frac{\partial a_{i}}{\partial z_j}|_{z=0}\right) c_i =0
\;\; \mbox{ for each } \;\; 1 \leq j \leq n.$$ It is straight
forward to check that   the point
$$z_1=\cdots = z_n =0, \;\;\;  [y_1:\cdots: y_{n-1}] = [c_1: \cdots :
c_{n-1}]$$ is a singular point  of the hypersurface $$ a_1(z_1,
\ldots, z_n)y_1 + \cdots + a_{n-1}(z_1, \ldots, z_n)
 y_{n-1}=0. \;\; \; \Box$$

\medskip
{\bf Corollary 4.6} {\it The irreducible hypersurface
$h^{-1}({\cal D})$ is the closure of the union of all rational
curves in $T^*(\S)$.}

\medskip
{\it Proof}. By Corollary 2.2 and Theorem 4.4 , $h^{-1}({\cal D})$
is covered by rational curves. By Proposition 4.2, there exists no
rational curve in $T^*(\S) \setminus h^{-1}({\cal D})$. $\Box$

\medskip
 {\bf Remark 4.7} Note that for a general
Hecke curve $\ell$,
$$ T(\S)|_{\ell} \cong {\cal O}(2) \oplus {\cal O}(1)^{2r-2} \oplus
{\cal O}^{(r^2-1)(g-1) - 2r +1}$$ because $\ell$ is a minimal
rational curve of $\S$ and $\ell \cdot K^{-1}_{\S} = 2r$. A
section $\tilde{\ell}$ of $T^*(\S)|_{\ell}$ gives a smooth
rational curve in $\hat{\cal S}$ in the notation of Corollary 2.2.
On the other hand, for a general $s \in {\cal D}$, $h^{-1}(s)$ is
an open subset of the compactified Jacobian of the spectral curve
of $s$ which has a unique node (cf. [KP, Remark 1.7]). The
normalization of the compactified Jacobian is a ${\bf P}_1$-bundle
over an abelian variety which is the Jacobian of the normalization
of the spectral curve. The curve $\tilde{\ell}$ is the image of a
fiber of this ${\bf P}_1$-bundle.  In [KP, Remark 1.7], it was
stated that the image of this ${\bf P}_1$-fiber has a node in the
compactified Jacobian. This is inaccurate, as $\tilde{\ell}$ is a
smooth rational curve. In fact, the morphism [KP, (1.13)] does not
exist because the pull-back of a torsion-free sheaf on the nodal
curve to its normalization is not necessarily torsion-free.  The
normalization map of the compactified Jacobian  identifies two
points on different fibers of the ${\bf P}_1$-bundle, contrary to
the claim in [KP, Remark 1.7]. However this mistake in Remark 1.7
does not affect the rest of the argument in [KP].

\medskip
{\bf Remark 4.8} It is possible to describe the section
$\tilde{\ell}$ in Remark 4.7 more explicitly as a family of Hecke
transformations of a Higgs field. Such description may give a more
direct proof of Corollary 4.6 and more detailed information about
the geometry of the Hitchin fibers. We will leave it for a future
investigation.

\medskip
{\bf Remark 4.9} In the above manner, the sections of $T^*(\S)$
over Hecke curves give a rank-1 foliation on an open subset of
$h^{-1}({\cal D})$. This foliation can be described in another
way. The cotangent bundle $T^*(\S)$ has a natural symplectic form.
The restriction of the symplectic form on the hypersurface
$h^{-1}(D)$ must be a holomorphic 2-form with a 1-dimensional
kernel. It is not difficult to check that the foliation defined by
this 1-dimensional kernel is precisely the foliation given by
Hecke curves.

\section{Applications}

As an application of Theorem 4.4, we will give a proof of the
non-abelian Torelli theorem, simplifying the proof in [KP, Theorem
E] for $g \geq 4$.

\medskip
 {\bf Theorem 5.1} {\it Let $C$ and $C'$ be two smooth
projective curves of genus $g \geq 4$. Let $f:{\cal SU}^s_C(r, d)
\rightarrow {\cal SU}^s_{C'}(r,d)$ be a biregular morphism. Then
$f$  induces a biregular morphism $C \cong C'$.}

\medskip
 The following is a direct consequence of Corollary 4.6.

\medskip
{\bf Lemma 5.2} {\it In the situation of Theorem 5.1, let
$$W:= H^0(C, \omega_C^{\otimes 2}) \oplus \cdots \oplus H^0(C,
  \omega_C^{\otimes r})$$
  $$W':= H^0(C', \omega_{C'}^{\otimes 2}) \oplus \cdots \oplus H^0(C',
  \omega_{C'}^{\otimes r})$$ be the spaces of characteristic polynomials and $$h:
T^*({\cal SU}^s_C(r,d)) \rightarrow W$$ $$h': T^*({\cal
SU}^s_{C'}(r,d)) \rightarrow W'$$ be the Hitchin maps.   Let
${\cal D} \subset W$ (resp. ${\cal D}' \subset W'$) be the
hypersurface consisting of characteristic polynomials with
singular spectral curves and ${\cal S} \subset {\bf P}T^*(\S)$
(resp. ${\cal S}' \subset {\bf P}T^*({\cal SU}^s_{C'}(r,d))$) be
the hypersurface corresponding to $h^{-1}({\cal D})$ (resp.
$(h')^{-1}({\cal D}')$). Let $df^* : {\bf P} T^*({\cal
SU}^s_{C'}(r,d)) \rightarrow {\bf P} T^*(\S) $ be the pull-back by
$f$. Then $ df^*({\cal S}') = {\cal S}$.}

\medskip
{\it Proof of Theorem 5.1}. Let $E$ be a general point of $\S$ and
$E'=f(E)$. By Lemma 5.2, $df^*_E({\cal S}'_{E'}) = {\cal S}_E$.
Thus by Theorem 4.4, $f$ induces a biregular morphism ${\cal C}_E
\cong {\cal C}_{E'}$. This induces a biregular morphism ${\bf P}
\Omega_E \cong {\bf P} \Omega_{E'}$ by Theorem 3.1 and Theorem
3.7, and consequently a biregular morphism $C \cong C'$ because
$C$ (resp. $C'$) is the Albanese image of ${\bf P}\Omega_E$ (resp.
${\bf P}\Omega_{E'}$). $\Box$

\medskip
A precise description of the automorphism group of $\S$, for $g
\geq 3$, was given by Kouvidakis and Pantev. An essential part of
their work was the following, which we will prove as an
application of Theorem 4.4. This simplifies the proof in [KP] for
$g \geq 4$.

\medskip
{\bf Theorem 5.3} {\it Let $C$ be a smooth projective curve of
genus $\geq 4$.  The group of automorphisms of $\S$ is generated
by automorphisms of the following two types when $r \not|\; 2d$.

(a) $E \mapsto \gamma^*E$ where $\gamma$ is an automorphism of the
curve $C$, and

(b) $E \mapsto E \otimes \mu$ where $\mu$ is an $r$-torsion of the
Picard group of $C$.

\smallskip When $r \; | \; 2d$,  additional generators of the following
type are needed.

(c)  $E \mapsto E^* \otimes \nu$ where $\nu$ is a line bundle of
degree $\frac{2d}{r}$ on $C$ whose $r$-th power is isomorphic to
the square of $\det(E)$.}

\medskip
We need two simple lemmas.

\medskip
{\bf Lemma 5.4} {\it Let $E$ (resp. $E'$) be a vector bundle of
rank $r$  on a smooth projective curve $C$  of genus $\geq 4$ and
$\Omega_E$ (resp. $\Omega_{E'}$) be the relative cotangent bundle
on ${\bf P}E^*$ (resp. ${\bf P}(E')^*$) with respect to the
natural projection $\varpi: {\bf P}E^* \rightarrow C$ (resp.
$\varpi': {\bf P} (E')^* \rightarrow C$). Suppose there exists a
biregular morphism $G: {\bf P} \Omega_E \rightarrow {\bf P}
\Omega_{E'}$. Then there exists a biregular automorphism $\gamma:
C \rightarrow C$ making the following diagram commutative. $$
\begin{array}{ccc} {\bf P} \Omega_E   &
\stackrel{G}{\longrightarrow} & {\bf P} \Omega_{E'} \\ \!
\downarrow
{\footnotesize \pi} & & \downarrow \! {\footnotesize \pi'} \\
C & \stackrel{\gamma}{\longrightarrow} & C
\end{array}$$
 Moreover, either $G$ descends to a biregular morphism $
{\bf P}E^* \rightarrow {\bf P} (E')^*$ or it descends to a
biregular morphism $ {\bf P}E^* \rightarrow {\bf P}E'$.}

\medskip
{\it Proof}. The existence of $\gamma$ is obvious by considering
Albanese map. Each fiber of $\pi$ and $\pi'$ is isomorphic to
${\bf P} T^*({\bf P}_{r-1})$ which has exactly two Mori
contractions (of extremal rays) $${\bf P}T^*({\bf P}_{r-1})
\longrightarrow {\bf P}_{r-1} \;\; \mbox{ and } \;\; {\bf
P}T^*({\bf P}_{r-1}) \longrightarrow {\bf P}_{r-1}^*.$$ Thus ${\bf
P} \Omega_E$ (resp. ${\bf P} \Omega_{E'}$) has exactly two Mori
contractions
$$ {\bf P} \Omega_E \longrightarrow {\bf P} E^* \;\; \mbox{ and }
\;\; {\bf P} \Omega_E \longrightarrow {\bf P} E $$
$$  \left( \mbox{ resp.} \;\; {\bf P} \Omega_{E'} \longrightarrow {\bf P}(E')^* \;\; \mbox{ and }
\;\; {\bf P} \Omega_{E'} \longrightarrow {\bf P}E'. \right) $$
Thus $G$ induces either ${\bf P}E^* \cong {\bf P} (E')^*$ or ${\bf
P}E^* \cong {\bf P} E'$. $\Box$

\medskip
{\bf Lemma 5.5} {\it In the situation of Lemma 5.4, assume that
$\deg(E) = \deg(E') =: d$ and $\det(E)= \gamma^* \det(E')$. Then
denoting by $\mbox{\rm Pic}^0(C)[r]$ the $r$-torsion subgroup of
the Picard group, one of the following holds.

(i) If $r \not| \;  2d$, there exists $\mu \in \mbox{\rm
Pic}^0(C)[r]$ such that $E \cong \gamma^*(E'  \otimes \mu)$.

(ii) If $r \; | \; 2d$, either there exists $\mu \in \mbox{\rm
Pic}^0(C)[r]$ such that $E \cong \gamma^*(E'  \otimes \mu)$, or
there exists $\nu \in \mbox{\rm Pic}^{\frac{2d}{r}}(C)$ with
$\nu^{\otimes r} = (\det(E'))^{\otimes 2}$ such that $E \cong
\gamma^*( (E')^* \otimes \nu).$}

\medskip
{\it Proof}. {}From Lemma 5.4, it is obvious that either $E \cong
\gamma^*( E' \otimes \mu)$ or $E \cong \gamma^*((E')^*  \otimes
\nu)$ for some line bundles $\mu, \nu$ on $C$. The assumption
$\det(E) = \gamma^*\det(E')$ can be easily translated to the
properties of $\mu$ and $\nu$ described in (i) and (ii). $\Box$

\medskip
{\it Proof of Theorem 5.3}.  Let $\sigma $ be an automorphism of
$\S$. Arguing as in the proof of Theorem 5.1,  we see that
$\sigma$ induces a biregular morphism $G: {\bf P} \Omega_E \cong
{\bf P} \Omega_{E'}$ for a general $E \in \S$ and $E'= \sigma(E)$.
By Lemma 5.4, $\sigma$ induces an automorphism $\gamma_E \in
Aut(C)$ for each general $E \in \S$. Since $Aut(C)$ is finite,
$\gamma_E$ is independent of $E$ for general $E$. Composing
$\sigma$ with an automorphism of type (a), we may assume that
$\gamma_E = I_C$, the identity map of $C$. Then $\sigma$ must
agree with an automorphism of type (b) or (c) by Lemma 5.5. $\Box$

\medskip
As a final application of Theorem 4.4, we prove the following
result on  morphisms between moduli spaces of bundles.

\medskip
 {\bf Theorem 5.6} {\it Let $C$ and $C'$ be two smooth projective
 curves of genus $g \geq 4$. Let $f: \S \rightarrow {\cal SU}^s_{C'}(r,d)$ be a
 surjective morphism. Then $f$ is biregular. }

\medskip
{\it Proof}.  The key point is to prove an analogue of Lemma 5.2.
In other words, when $df^*: f^*T^*({\cal SU}^s_{C'}(r,d))
\rightarrow T^*(\S)$ is the natural morphism associated to $f$, we
claim that  $df^*(f^* {\cal S}') \subset {\cal S}$ in the notation
of Lemma 5.2. Note that the proof of Lemma 5.2 does not work when
$f$ is {\it a priori} not biregular.

To prove the claim, let $\ell$ be a general Hecke curve on ${\cal
SU}^s_{C'}(r,d)$ and $\hat{\ell} \subset \S$ be an irreducible
component of $f^{-1}(\ell)$. An element $\sigma \in H^0(\ell,
T^*({\cal SU}^s_{C'}(r,d))$ defines an element $f^*\sigma \in
H^0(\hat{\ell}, T^*(\S))$. Let $\ell^{\flat}$ be the image of
$f^*\sigma$ in $T^*(\S)$. Since $W$ is affine, $h(\ell^{\flat})$
is one point $s \in W$. Suppose $s \not\in {\cal D}$. By
Proposition 4.2, $h^{-1}(s)$ is an open subset $A'$ of an abelian
variety $A$. The natural projection $A' \rightarrow \S$ is
dominant and so is its composition with $f$, which is denoted by
$f': A' \rightarrow {\cal SU}^s_{C'}(r,d)$. But the complete curve
$\ell^{\flat}$ satisfies $f'(\ell^{\flat}) = \ell$. This is a
contradiction to Proposition 2.4, because the variety of Hecke
tangents is not linear. It follows that $s \in {\cal D}$. Thus
$f^* \sigma$ has its image in $h^{-1}({\cal D})$. Since
$(h')^{-1}({\cal D}')$ is covered by images of $\sigma$ by
Corollary 2.2 and Theorem 4.4, this implies that
$df^*((h')^{-1}({\cal D}')) \subset h^{-1}({\cal D})$, as claimed.

Choose an analytic open subset $U \subset \S$ such that   $f|_U: U
\rightarrow f(U)$ is biholomorphic. By the claim, for each $u \in
U, (df_u)^*({\cal S}'_{f(u)}) = {\cal S}_u$. By Theorem 4.4,
$df_u({\cal C}_u) = {\cal C}'_{f(u)}$. Then we can proceed as in
the proof of Theorem 5.1 and Theorem 5.3 to show that $C \cong C'$
and $f|_U$ agrees with the restriction of an automorphism of $\S$
to $U$. Hence $f$ is biregular. $\Box$.

\medskip
{\bf Remark 5.7} The argument used in the proof of Theorem 5.6
shows that when $s \not\in {\cal D}$, the projection $h^{-1}(s)
\rightarrow \S$ cannot be proper over a Hecke curve. In fact,
there exists no complete curve in $h^{-1}(s)$ which is mapped to a
Hecke curve in $\S$. Since for any subvariety $Z$ of codimension
$\geq 2$ in $\S$, there exists a Hecke curve disjoint from $Z$,
this means that the locus where  the projection $h^{-1}(s)
\rightarrow \S$ is not proper is of codimension 1 in $\S$.

\bigskip
{\bf Acknowledgment} The first named author would like to thank
Dr. X. Sun for a discussion regarding Remark 4.7.  The second
named author wishes to thank KIAS for the hospitality extended to
him in Spring 2003, during which this collaboration took place.

\bigskip
{\bf References}

\medskip

[An] Andreotti, A.: On a theorem of Torelli. Amer. J. Math. {\bf
80} (1958) 801--828

 [BNR] Beauville, A., Narasimhan, M.S. and Ramanan, S.:
Spectral curves and the generalized theta divisor. J. reine angew.
Math. {\bf 398} (1989) 169-179

[Hi] Hitchin, N. J.: Stable bundles and integrable systems. Duke
Math. J. {\bf 54} (1987) 91-114

[HM1] Hwang, J.-M. and Mok, N.: Projective manifolds dominated by
abelian varieties. Math. Z. {\bf 238} (2001) 89-100

[HM2] Hwang, J.-M. and Mok, N.: Finite morphisms onto Fano
manifolds of Picard number 1 which have rational curves with
trivial normal bundles. J. Alg. Geom. {\bf 12} (2003) 627-651

[HM3] Hwang, J.-M. and Mok, N.: Birationality of the tangent map
for minimal rational curves.  alg-geom/0304101, to appear in Asian
J. Math., {\it Special issue dedicated to Yum-Tong Siu on his 60th
birthday}

 [Hw1] Hwang,
J.-M.: Tangent vectors to Hecke curves on the moduli space of rank
2 bundles over an algebraic curve. Duke Math. J. {\bf 101} (2000)
179-187

[Hw2] Hwang, J.-M.: Geometry of minimal rational curves on Fano
manifolds. {\it School on Vanishing Theorems and Effective Results
in Algebraic Geometry (Trieste, 2000)}, 335--393, ICTP Lect.
Notes, 6, Abdus Salam Int. Cent. Theoret. Phys., Trieste, 2001

 [Hw3] Hwang, J.-M.: Hecke curves on the moduli space of vector
bundles over an algebraic curve. Proceedings of the symposium {\it
Algebraic Geometry in East Asia (Kyoto, August 3-10, 2001)},
155-164, World Scientific, 2003

 [KP] Kouvidakis, A. and Pantev, T.: The automorphism group of the
 moduli space of semi-stable bundles. Math. Annalen {\bf 302}
 (1995) 225-268

[La] Laumon, G.: Un analogue global du c\^one nilpotent. Duke
Math. J. {\bf 57} (1988) 647-671

 [NR1] Narasimhan, M.S. and Ramanan, S.: Deformations of the
moduli space of vector bundles over an algebraic curve. Ann. Math.
{\bf 101} (1975) 391-417

[NR2] Narasimhan, M.S. and Ramanan, S.: Geometry of Hecke cycles
I. in {\em C. P. Ramanujam-a tribute}. Springer Verlag, 1978,
p.291-345

\bigskip
\begin{tabbing}
aaaaaaaaaaaaaaaaaaaaaaaaaaaaaaaaaaaaaaaa\= \kill
Jun-Muk Hwang \> S. Ramanan \\

Korea Institute for Advanced Study  \> Institute of Mathematical
Sciences \\
207-43 Cheongryangri-dong \>  Taramani, Chennai 600, India $\;$ {\it and}\\

Seoul 130-722, Korea  \> Chennai Mathematical Institute\\

e-mail: jmhwang@ns.kias.re.kr \> G.N. Chetty Road, Chennai 600
017, India
\\ \>  e-mail: sramanan@cmi.ac.in
\end{tabbing}

\end{document}